\documentclass[12pt]{article}

\usepackage{amscd,amsmath, amssymb, fancyhdr,
	epsfig,color,bbm,url}
\definecolor{dblue}{rgb}{0,0,.6}

\usepackage[backref=page]{hyperref}
\usepackage[T1]{fontenc}
\renewcommand*{\backref}[1]{}
\renewcommand*{\backrefalt}[4]{%
	\ifcase #1 (Not cited.)%
	\or        (Cited on page~#2.)%
	\else      (Cited on pages~#2.)%
	\fi}

\hypersetup{
	colorlinks   = true,
	citecolor    = magenta,
	linkcolor= blue
}

\numberwithin{equation}{section}

\def\eqref#1{(\ref{#1})}

\newcommand{\Z}{{\mathbb Z}}

\newcommand{\R}{{\mathbb R}}
\newcommand{\Q}{{\mathbb Q}}
\renewcommand{\L}{\mathcal{L}}

\def\1{\sqrt{-1}\.}

\newcommand{\cntrct}                
{\hspace{2pt}\raisebox{1pt}{\text{$\lrcorner$}}\hspace{2pt}}

\makeatletter

\renewcommand{\bar}{\overline}
\renewcommand{\phi}{\varphi}
\renewcommand{\epsilon}{\varepsilon}

\newcommand{\Id}{\operatorname{Id}}

\renewcommand{\Im}{{\operatorname{Im}}}

\renewcommand{\d}{\partial}

\newcommand{\ds}{\frac{\partial}{\partial s}}

\usepackage{amsthm}
\theoremstyle{definition}
\newtheorem{theorem}{Theorem}[section]
\newtheorem{lemma}[theorem]{Lemma}
\newtheorem{proposition}[theorem]{Proposition}

\newtheorem{cor}[theorem]{Corollary}
\newtheorem{defin}[theorem]{Definition}
\newtheorem{example}[theorem]{Example}
\theoremstyle{remark}
\newtheorem{rem}{Remark}[section]

\makeatletter

\@addtoreset{equation}{section} \@addtoreset{footnote}{section} \makeatother

\begin{document}
	\begin{center}
		{\LARGE\bf
			Locally conformally Hessian and statistical manifolds.}
		\medskip
		\medskip

		Pavel Osipov\footnote{National Research University Higher School of Economics, Russian Federation } \footnote{Pavel Osipov is 
			partially supported by the HSE University Basic
			Research Program and by the Theoretical Physics and Mathematics
			Advancement Foundation "BASIS" .}
	\end{center}
	
	%
	
		\begin{abstract}
			A statistical manifold $\left(M,D,g\right)$ is a manifold $M$ endowed with a torsion-free connection $D$ and a Riemannian metric $g$ such that the tensor $D g$ is totally symmetric. If $D$ is flat then $\left(M,g,D\right)$ is a Hessian manifold.  
			A locally conformally Hessian (l.c.H) manifold is a quotient of a Hessian manifold $(C,\nabla,g)$ such that the monodromy group acts on $C$ by Hessian homotheties, i.e. this action preserves $\nabla$ and multiplies $g$ by a group character. The l.c.H. rank is the rank of the image of this character considered as a function from the monodromy group to real numbers. A l.c.H. manifold is called radiant if the Lee vector field $\xi$ is Killing and satisfies $\nabla \xi =\lambda \Id$.  We prove that the set of radiant l.c.H. metrics of l.c.H. rank 1 is dense in the set of all radiant l.c.H. metrics. We prove a structure theorem for compact radiant l.c.H. manifold of l.c.H. rank 1. Every such manifold $C$ is fibered over a circle, the fibers
		are statistical manifolds of constant curvature, the fibration is locally trivial, and $C$ is reconstructed from the statistical structure on the fibers and
		the monodromy automorphism induced by this fibration.
	\end{abstract}

		\tableofcontents
	\section{Introduction}
	 A {\bfseries flat affine manifold} is a differentiable manifold equipped with a flat torsion-free connection. Equivalently, it is a manifold equipped with an atlas such that all transition functions between charts are affine transformations (see \cite{FGH} or \cite{Sh}). 
	A {\bfseries Hessian manifold} is a flat affine manifold $\left(C,\nabla\right)$ with a metric $g$ which is locally equivalent to a Hessian of a function. Equivalently, the metric $g$ is Hessian if and only if the tensor $\nabla g$ is totally symmetric.

	
	The metric $g$ on a flat affine manifold is called {\bfseries locally conformally Hessian (l.c.H)} if for any open neighborhood $U\subset C$ there exists a function $f$ on $U$ such that the locally defined metric $e^{-f} g$ is Hessian. The main purpose of this paper is to describe compact flat affine manifolds with an l.c.H. metric.
	
	A {\bfseries statistical manifold} $\left(C,D,g\right)$ is a manifold $M$ endowed with a torsion-free connection $D$ and a Riemannian metric $g$ such that the tensor $D g$ is totally symmetric. The term statistical manifolds arose in information geometry (see \cite{AN}). In this sense, statistical manifolds is a space of probability distributions endowed with the Fisher information metric. For example, the statistical manifold corresponding to the family of normal distributions is isometric to the hyperbolic plane. 
	
	 A statistical manifold $\left(C,D,g\right)$ is said to be {\bfseries of constant curvature} $c$ if the curvature tensor $\Theta_D$ satisfies 
	$$
	\Theta_D (X,Y)Z=c\left(g(Y,Z) X - g(X,Z) Y\right),
	$$ 
	for any $X,Y,Z\in TM$. For example, a Riemannian manifold of constant sectional curvature is a statistical manifold of constant curvature. The definition of statistical manifolds of constant curvature arose in the context of geometry of affine hypersurfaces (\cite{K}). We describe this origin in Section \ref{ku} (Theorems \ref{K0} and \ref{k}). Note that Hessian manifolds are statistical manifolds of curvature 0. We assume that the curvature of a statistical manifold is not equal to 0.

	Convex projective geometry provides a wide class of statistical manifolds. A domain $U\subset \mathbb{RP}^n$ is called properly convex if the closure of $U$ is a compact convex set in some affine chart. If $\Gamma$ is a discrete subgroup of the group of projective automorphisms of a properly convex domain $U\subset \mathbb{RP}^n$ such that $M=U/\Gamma$ is a manifold then $M$ is called a {\bfseries properly convex} $\mathbb{RP}^n$-manifold. For examples of compact properly convex $\mathbb{RP}^n$-manifolds see \cite{B}.
	\begin{theorem}[\cite{KO}]\label{OT}
		Any properly convex $\mathbb{RP}^n$-manifold admits a statistical structure of negative constant curvature. Any compact statistical manifold of negative constant curvature admits a properly convex $\mathbb{RP}^n$ structure.
	\end{theorem} 
The proof of Theorem \ref{OT} is based on the characteristic of properly convex $\mathbb{RP}^n$-structures on compact manifolds given in \cite{L}. This characteristic is equivalent to the following theorem. 
\begin{theorem}[\cite{L} and \cite{KO}]
	A compact statistical manifold $(M,D,g)$ of negative constant curvature admits a properly convex $\mathbb{RP}^n$ structure if and only if $M$ admits a $D$-flat volume form.
\end{theorem}
	We will give an alternative proof of Theorem \ref{OT} using the results of the present paper.  
	
		A {\bfseries self-similar manifold} $(M,g,\xi)$ is a Riemannian manifold endowed with a vector field $\xi$ satisfying $\L_\xi g =2g$ (\cite{A}). A {\bfseries self-similar Hessian manifold} $\left(C,\nabla,g,\xi\right)$ is a Hessian manifold $(C,\nabla,g)$ endowed with an affine vector field $\xi$ such that $(C,g,\xi)$ is a self-similar manifold (\cite{Os}, Definition 3.3). 
		If the field $\xi$ satisfies $\nabla \xi = \lambda \text{Id}$ then we say that $\left(M,\nabla,g,\xi\right)$ is a {\bfseries radiant Hessian manifold} ((\cite{Os}, Definition 3.8)). The condition $\xi = \lambda \text{Id}$ means that there exists a flat affine atlas on $C$ such that in local coordinates we have $\xi=\lambda\sum \frac{\d}{\d x^i}$, i.e. $\xi$ is proportional to the radiant vector field. The field $\xi$ on a self-similar Hessian manifold $\left(C,\nabla,g,\xi\right)$ is a potential vector field (i.e. $\xi$ equals the gradient of a function) if and only if $\left(C,\nabla,g,\xi\right)$ is a direct product of radiant Hessian manifolds (\cite{Os}, Theorem 1.4).
		
		The total space of the tangent bundle of a Hessian manifold can be endowed with a K\"ahler structure (\cite{Sh}). A globally conformally K\"ahler structure on the total space of the tangent bundle of a selfsimilar Hessian manifold was constructed in \cite{Os2}.

There is a correspondence between radiant Hessian manifolds and statistical manifolds of constant curvature. Precisely, a Riemannian cone $$\left(M\times \R^{>0},s^2g_M+ds^2\right)$$ over a statistical manifold $\left(M,g_M,D\right)$ of constant curvature admits a structure of a radiant Hessian manifold. Conversely, level sets of a Hessian potential on a radiant Hessian manifold are statistical manifolds of constant curvature (\cite{Os3}). 
	
	By $dd^c$ Lemma, any K\"ahler form can be locally represented as a complex Hessian $dd^c \varphi$. Hence, we can consider Hessian manifold as a real analogue of K\"ahler manifolds. 	A {\bfseries Sasakian manifold} is a Riemannian manifold $(M,g_M)$ such that the cone metric ${g=s^2 g_M+ds^2}$  on $M\times\R^{>0}$ is K\"ahler with respect to a dilatation-invariant complex structure $I$ (see \cite{OV2} or \cite{OV4}). Thus, we can consider statistical manifolds of constant curvature as an analogue of Sasakian manifolds. 
	
	We extend this analogy and define {\bfseries locally conformally Hessian (l.c.H.) manifolds} similarly to locally conformally K\"ahler manifolds (see e.g. \cite{OV3}). An l.c.H. manifold $(C,\nabla, g,\theta)$ is a flat affine manifold $(C,\nabla)$ endowed with a Riemannian metric $g$ and closed 1-form $\theta$ such that $\nabla g-\theta\otimes g$ is a totally symmetric tensor. The form $\theta$ and the vector field $\xi =\theta^\sharp$ are called a {\bfseries Lee form} and a {\bfseries Lee vector field}. If locally $\theta = df$ then $e^{-f} g$ is a locally defined Hessian metric.
	
	We study l.c.H. manifolds with an affine Killing Lee vector field. This class of manifolds is an analogue of Vaisman manifolds: a l.c.K. manifold is called {\bfseries Vaisman} if the Lee vector field is Killing (\cite{Va1}. The Lee vector field of any Vaisman manifold is holomorphic. The affine structure on an l.c.H. manifold take the same place as the complex structure on a l.c.K manifold. In contrast to the l.c.K. case, the Lee vector field on an l.c.H. manifold can be Killing but not affine (Example \ref{e67}).
	
	Note that Vaisman manifolds belong to the larger class of l.c.K. manifolds called l.c.K. manifolds with a potential. A l.c.K. manifold is called a l.c.K. manifolds with a potential if it admits a Kähler covering on which the Kähler metric has global, positive and proper potential function (\cite{OV3},\cite{OV5},\cite{OV6}).

	Let $\left(C,\nabla,g,\theta\right)$ be a l.c.H. manifold and $\xi$ its Lee vector field. Then we say that $\left(C,\nabla,g,\theta\right)$ is a {\bfseries radiant l.c.H. manifold} if $\xi$ is Killing and there exists a constant $\mu\in \R$ such that $\nabla \xi=\mu\Id$. Equivalently, a l.c.H. manifold is radiant if it is a quotient of a radiant Hessian manifold $\left(\widetilde C, \widetilde \nabla, \widetilde g\right)$ such that the monodromy group acts on $\widetilde C$ by Hessian homotheties (see Definition \ref{414} and Proposition \ref{415}). The following theorem motivates us to consider radiant l.c.H. manifolds.  

\begin{theorem}\label{TT3}
	Let $\left(C,\nabla,g,\theta\right)$ be a compact l.c.H. manifold with an affine Killing Lee vector field $\xi$. Suppose $({C},{g})$ is not a locally conformally flat Riemannian manifold. Then $\left(C,\nabla,\theta,\xi\right)$ is a radiant Hessian manifold.
\end{theorem} 
	
	
	

	An l.c.H. manifold $(C,\nabla,g,\theta)$ admits a covering $\widetilde{C}$ endowed with a Hessian metric $\widetilde{g}$ such that the deck group $\text{Aut}_C(\widetilde{C})$ acts on $\widetilde{C}$ by Hessian homotheties. 
	It defines a character 
	$$
	\chi : \text{Aut}_C(\widetilde{C}) \to \R^{>0}.
	$$
	Similarly to Definition 2.11 in \cite{GOPP}, we call the rank of the group ${\Im\chi\subset \R^{>0}}$ by {\bfseries l.c.H. rank} of  $(C,\nabla,g,\theta)$. The l.c.H. rank of an l.c.H. manifold $(M,\nabla,g,\theta)$ is equal to 1 if and only if $[\theta]\in H^1(M,\mathbb Q)$. 
	
	Now we describe the main results of the present paper.
	
	\begin{theorem}\label{14}
		Let $\left(C,\nabla\right)$ be a compact flat affine manifolds and ${{L}\subset H^1(C,\R)}$ a set of cohomology classes of Lee forms of radiant l.c.H. structures. Then $L$ is open.  
	\end{theorem}

	It follows from Theorem \ref{14} that any radiant l.c.H. structure can be approximated by a radiant l.c.H. structure of rank 1. 
	
	\begin{cor}\label{ccc}
		Let $(C,\nabla,g,\theta)$ be a radiant l.c.H. manifold. Then there exist a metric $g'$ and a 1-form $\theta'$ on $C$ such that $(C,\nabla,g',\theta')$ is a radiant l.c.H. manifold of l.c.H. rank 1.
	\end{cor}

	\begin{theorem}\label{t1}
		Let $\varphi$ be an automorphism of a statistical manifold $(M,g_M,D)$  of constant curvature. Consider the automorphism 
		$$
		{\varphi_q: M\times \R^{>0}\to M\times \R^{>0}, \ \ \ \varphi_q(m,t)= \left(\phi(m),qt\right)}.
		$$
		Then $M\times \R^{>0} /_{\varphi_q}$ admits a radiant l.c.H. structure of l.c.H. rank 1. 
		\end{theorem}
	The construction from Theorem \ref{t1} is similar to the construction of a Vaisman manifold by a Sasakian manifold described in \cite{OV1}
	
	If the Lee vector field $\xi$ on a l.c.H. manifold $(C,\nabla,g,\xi)$ is Killing then $a:=g(\xi,\xi)$ is a constant (Corollary \ref{69}). 
	\begin{theorem}\label{t2}
		Let $(C,\nabla,g,\theta)$ be  a compact radiant l.c.H. manifold of l.c.H. rank 1, $\xi$ the Lee vector field, and $\nabla \xi =\mu \Id$, where $\mu$ is a constant. Then  $(C,\nabla,g,\theta)$ can be constructed from a statistical manifold $\left(M,D,g_M\right)$ of constant curvature as in Theorem \ref{t1}. Moreover, $(M,D,g)$ is a statistical manifold of negative constant curvature if and only if $\mu\in (-\infty, -a)\cup (0,\infty)$, where $a=g(\xi,\xi)$. 
	\end{theorem}
	Theorems \ref{14}, \ref{t1}, and \ref{t2} are analogues to the structure theorems for compact Vaisman manifolds from \cite{OV3} (except the condition on sign of curvature).
	
	Combining Theorem \ref{OT}, Corollary \ref{ccc}, and Theorem \ref{t2} we get the following corollary.

	\begin{theorem}\label{t4}
		Let $(C,\nabla,g,\theta)$ be a radiant l.c.H. manifold, $\xi$ the Lee vector field, and $\nabla \xi =\mu \Id$, where $\mu$ is a constant. Suppose ${}\mu\in (-\infty -a)\cup (0,\infty)$, where $a=g(\xi,\xi)$. Then the universal covering of $C$ is isomorphic (as a flat affine manifold) to a convex cone without full straight lines.
	\end{theorem}

\section{Self-similar manifolds}
\subsection{Self-similar manifolds: definition and basis properties.}

\begin{defin}
	A {\bfseries self-similar manifold} $(C,g,\xi)$ is a Riemannian manifold $(C,g)$ endowed with a field $\xi$ satisfying  
	$$
	\L_\xi  g = 2g.
	$$
	If $\xi$ is complete then the manifold is called a {\bfseries global self-similar manifold}.
\end{defin}

It follows from the definition that a global self-similar manifold is a Riemannian manifold endowed with a 1-parameter group of homothetic automorphims $\{\varphi_t\}$ such that $\varphi_t^* g= e^{2t} g$. The term "self-similar" is motivated by the fact that for any $\lambda\in \R^{>0}$, a global self-similar manifold $(C,g)$ is isometric to $(C,\lambda g)$. 

\begin{example}
	Let $\left(C=M\times\R^{>0},g=s^2g_M+ ds^2\right)$ be a {\bfseries Riemannian cone} and $\xi=s\frac{\d}{\d s}$. Then $\left(C,g,\xi\right)$ is a global self-similar manifold.
\end{example}
\begin{example}[\cite{Os}]
	Let $\varphi$ and $s$ be coordinates on $S^1$ and $\R^{>0}$. Then the collection $\left(C=S^1\times \R^{>0},g=s^2d\varphi^2+sds\cdot d\varphi+ds^2,s\ds\right)$ is a global self-similar manifold but  $\left(C,g\right)$ is not isometric to a Riemannian cone.
\end{example}

\begin{defin}
	We say that $(C,g,\xi)$ is a {\bfseries self-similar manifold with a potential homothetic vector field} if $(M,g,\xi)$ is a self-similar manifold and $\xi$ is locally defined as a gradient of a function. If $\xi=\text{grad} \ f$ on a domain $U$ then $\iota_\xi g|_U=df$. Moreover, a form is closed if and only if it is locally exact. Therefore, the vector field $\xi$ is potential if and only if $d\iota_\xi g=0.$
\end{defin}

\begin{theorem}[\cite{Os}]\label{T1}\ 
	Let $\left(C,g,\xi\right)$ be a global self-similar manifold with a potential homothetic vector field.
	\begin{itemize} 
		\item [(i)] If $\xi$ vanishes at a point then $\left(C,g,\xi\right)$ is a Euclidean space with a radiant vector field $\left(\R^n,\sum_{i=1}^n\left(dx^i\right)^2,\sum x^i\frac{\d}{\d x^i} \right)$.  
		\item [(ii)] If $\xi$ does not vanish at any point then $\left(C,g,\xi\right)$ is a Riemannian cone $\left( M\times \R^{>0},s^2g_M+ds^2,s\ds \right)$. 
	\end{itemize}
\end{theorem}

		\subsection{Self-similar Hessian manifolds}\label{s3}
	
	\begin{defin}\label{d}
		A {\bfseries flat affine manifold} $(C,\nabla)$ is a differentiable manifold $C$ equipped with a flat torsion-free connection $\nabla$. Equivalently, it is a manifold equipped with an atlas such that all transition maps between charts are affine transformations (see e.g. \cite{FGH}). A {\bfseries radiant manifold} $(C,\nabla, \rho)$ is a flat affine manifold $(C,\nabla)$ endowed with a {\bfseries radiant vector field } $\rho$ i.e. a vector field satisfying
		\begin{equation}\label{3.1}
		\nabla \rho =\text{Id}
		\end{equation}
		Equivalently, it is a manifold equipped with an atlas such that all transition maps between charts are linear transformations. In the corresponding local coordinates we have 
		$$
		\rho=\sum x^i \frac{\d}{\d x^i}
		$$ 
		(see e.g. \cite{Go}).  
	\end{defin}

\begin{defin}{\label{Hes}}
	A {\bfseries Hessian manifold} $(C,\nabla, g)$ is a flat affine manifold endowed with a Riemannian metric such that $\nabla g$ is a totally symmetric tensor.
\end{defin}

	\begin{defin}
	Let $(C,\nabla)$ be a flat affine manifold. A vector field $\xi$ on $C$ is called {\bfseries affine} if the flow along $\xi$ preserves the connection $\nabla$.
\end{defin}	
\begin{proposition}[\cite{Go}]
	Let $X$ be a vector field on a flat affine manifold $(C,\nabla)$ then the following conditions are equivalent:
	\begin{itemize}
		\item[(i)] $X$ is affine.
		\item[(ii)] $\nabla X$ is a $\nabla$-flat 1-1 tensor.
		\item[(iii)] For any local coordinates $(x^1,\ldots, x^n)$, we have $X=\sum a_i(x) \frac{\d}{\d x^i}$, where $a_1(x),\ldots,a_n(x)$ are linear functions.
	\end{itemize}
\end{proposition}
	
	\begin{defin}
		A {\bfseries self-similar Hessian manifold} $(C,\nabla, g, \xi)$ is a Hessian manifold $(C,\nabla,g)$ endowed with an affine vector field $\xi$ such that $\L_\xi g =2g$ and the flow along $\xi$ preserves $\nabla$. 
	\end{defin} 


	\begin{defin}[\cite{Os}] \label{ddd}
		A self-similar Hessian manifold $(C,\nabla, g, \xi)$ is a {\bfseries radiant Hessian manifold} if it admits a radiant vector field $\rho$ and a constant $\lambda\in \R $ such that $\xi=\lambda\rho$. 
		Equivalently, there is a flat affine atlas such that in the corresponding local coordinates we have 
		$$
		\xi=\lambda\sum  x^i \frac{\d}{\d x^i}
		$$ 
		(see Definition \ref{d}). 
	\end{defin}

\begin{rem}\label{rem}
	The case $\lambda=2$ is studied in \cite{G-A} in the context of Equilibrium Thermodynamics. In this case, we have $\L_\rho g=g$ where $\rho= \sum  x^i \frac{\d}{\d x^i}$ and $g=\text{Hess}\ \varphi$ locally. Then $\varphi$ is linear along $\rho$ and $\iota_\rho g =\iota_\rho (\text{Hess}\ \varphi)=0$. Hence, if $\lambda=2$ then $g$ cannot be positive definite. If $\lambda =0$ then $\L_\xi g=0\ne 2g$. Thus, $\lambda \in \R \setminus \{0,2\}$.
\end{rem}

	We say that a self-similar Hessian manifold $\left(U,\nabla,\xi,g\right)$ is a direct product of self-similar Hessian manifolds $\left(U_i,\nabla_i,\xi_i,g_i\right)$ if
	$$U=\prod U_i, \ \ \ \nabla=\sum \nabla_i, \ \ \ g= \sum g_i, \ \ \ \xi=\sum \xi_i. $$
	Since the condition $\forall i\in \{1\ldots n\}: \L_{\xi_i} g_i =2g_i$ implies $\L_\xi g=2g$, the product of self-similar Hessian manifold is a self-similar Hessian manifold.

\begin{theorem} \label{t01}
	Let $\left(C,\nabla,g,\xi\right)$ be a self-similar Hessian manifold. Then the following conditions are equivalent.
	\begin{itemize}
		\item[(i)] The vector field $\xi$ is potential.
		\item[(ii)] There are $\lambda_1,\ldots,\lambda_k\in \R$ and a $\nabla$-flat decomposition of the tangent bundle $TC=\bigoplus_{i=1^k} V_i$ such that for any $i\in\{1,\ldots,k\}$ we have ${\nabla \xi|_{V_i}=\lambda_i \text{Id}}.$ 
	\item[(iii)] For any point $p\in C$ there is a neighborhood $U\ni p$ and a collection of radiant Hessian manifolds $\left(U_i,\nabla_i,\xi_i,g_i\right)$ such that $\left(U,\nabla,\xi,g\right)$ is a direct product of self-similar Hessian manifolds $\left(U_i,\nabla_i,\xi_i,g_i\right)$ and the decomposition to the direct product $U=\prod U_i$ corresponds to ${TU=\oplus V_i|_U}$. 
	\end{itemize}
\end{theorem}
 Theorem \ref{t01} can be reformulated in a short form.
	\begin{theorem} \label{t00}
	Let $\left(C,\nabla,g,\xi\right)$ be a self-similar Hessian manifold. Then $\xi$ is potential if and only if $\left(C,\nabla,g,\xi\right)$ is locally isomorphic to a direct product of radiant Hessian manifolds.
\end{theorem}

\begin{cor}\label{c1}
	Let $\left(C,\nabla,g,\xi\right)$ be a self-similar Hessian manifold with a potential homothetic vector field $\xi$. Suppose the holonomy of $(C,g)$ is irreducible. Then $\left(C,\nabla,g,\xi\right)$ is a radiant Hessian manifold.   
\end{cor}

	\section{Statistical manifolds}
	\subsection{Statistical manifolds: definition and examples}\label{s4}

	\begin{defin}
		A {\bfseries statistical manifold} $\left(M,D,g\right)$ is a manifold $M$ endowed with a torsion-free connection $D$ and a Riemannian metric $g$ such that the tensor $D g$ is totally symmetric. A statistical manifold $\left(M,D,g\right)$ is said to be {\bfseries of constant curvature} $c\in \R$ if the curvature tensor $\Theta_D$ satisfies 
		\begin{equation}\label{e}
		\Theta_D (X,Y)Z=c\left(g(Y,Z) X - g(X,Z) Y\right),
		\end{equation}
		for any $X,Y,Z\in TM$ (\cite{K}).
	\end{defin}
The following example motivates the definition above.
\begin{example}\label{ex1}
	Let $D$ be the Levi-Civita connection on a Riemannian manifold $\left(M,g\right)$. Then $\left(M,g,D\right)$ is a statistical manifold. The sectional curvature of $M$ is constant and equals $c\in \R$ if and only Equation \eqref{e}  is satisfied. Thus, $\left(M,g,D\right)$ is a statistical manifold of constant curvature if and only if the sectional curvature of $(M,g)$ is constant.  
\end{example}
\begin{example} \label{e33}
	Let $V\subset \R^{n+1}$ be a convex cone without full straight lines and 
	$$
	V^*=\{y\in \left(\R^{n+1}\right)^* | \ \forall x\in V: (x, y)>0\}
	$$ 
	be the dual cone. Consider the {\bfseries characteristic function}
	$$
	\psi(x) = \int_{V^*} e^{-(x,y)} dy
	$$
	and the {\bfseries characteristic hypersurface} $S=\{x\in V| \psi(x)=1\}$. Then the bilinear function $\psi$ is smooth and $g=\text{Hess} \left(\ln \psi\right)$ is positive definite (\cite{V}). Let $\nabla$ be a standard connection on $\R^{n+1}$ and $\xi=\sum_{i=1}^{n+1}x^i\frac{\d}{\d x^i}$. Set a connection $D$ on $S$ and a bilinear form $h$ on $S$ as the projection of $\nabla$ on components of the decomposition $TV|_S=TS\oplus \R\xi$, i.e. $D$ and $h$ are the connection $D$ on $S$ and the bilinear form $h$ on $S$ satisfying  
	$$
	\forall X,Y \in TS \ \ : \ \ \nabla(X,Y) =D_X Y +h(X,Y) \xi.
	$$   
	Then $\left(S,g|_S,D\right)$ is a statistical manifold of constant negative curvature (\cite{Sh}, Example 5.1 and Corollary 5.3). Let $\text{Aut}_\text{SL}(V)$ be the subgroup of automorphisms of $V$ in $\text{SL}\left(\R^{n+1}\right)$. Then $\text{Aut}_\text{SL}(V)$ preserves $\left(S,g|_S,D\right)$ (\cite{V}). 
\end{example}
\begin{rem}
	A convex cone $V$ is called homogeneous if there is a transitive action on $V$ by a group of linear automorphisms. The characteristic function of a homogeneous cone is a solution of a Monge-Ampère equation and characteristic hypersufaces of cones are affine spheres of a negative constant curvature (\cite{S}).  
\end{rem}
The following example is a reformulation of the previous one.
	\begin{example}\label{e34}
		We say that a domain $\Omega\subset \R \mathbb P^n$ is a convex domain if the intersection of $\omega$ with any full projective line is connected. Let $\Omega\subset \R \mathbb P^n$ be a convex domain without any full projective line, $\text{Aut}(\Omega)$ the group of projective automorphisms of $\omega$, and $\Gamma$ a discrete subgroup of $\text{Aut}(\Omega)$ such that $M=\Omega /\Gamma$ is a compact manifold. Then $M$ is called a {\bfseries properly convex $\R \mathbb{P}^n$ manifold} (see \cite{B}). Consider the projection $\pi:\R^{n+1}\setminus\{0\} \to \R\mathbb P^n$. A connected component $V\subset \R^{n+1}$ of the preimage $\pi^{-1}\Omega$ of a properly convex  domain $\Omega$ is a convex cone without straight full lines. We can identify $\Omega$ with a characteristic hypersurface $S$ of $V$ and $\Gamma$ with a subgroup of $\text{Aut}_\text{SL}(V)$. Hence, the $\text{Aut}_\text{SL}(V)$-invariant statistical structure on $S$ from Example \ref{e33} can be identified with a $\text{Aut}(\Omega)$-invariant statistical structure on $U$. Therefore, this statistical structure of constant curvature on $S\simeq U$ can be pushed forward to $M=U/\Gamma$. That is, any properly convex  $\R \mathbb{P}^n$ manifold admit a statistical structure of constant negative curvature.
\end{example}

It is proved in \cite{KO} that any compact statistical manifold of negative constant curvature  is a properly convex $\mathbb{RP}^n$ manifold.

\begin{example}[\cite{FHOSS}]
	A {\bfseries Sasakian manifold} is a Riemannian manifold $(M,g_M)$ such that the cone metric $ g=s^2 g_M+ds^2$  on $M\times\R^{>0}$ is K\"ahler with respect to a dilatation-invariant complex structure $I$ (see \cite{OV2} or \cite{OV4}). The field $\xi={I\ds \in T\left(M\times 1\right) \simeq TM}$ is called the {\bfseries Reeb vector field}. For any $X,Y\in TM$ set $K(X,Y)=g(X,\xi)g(Y,\xi)\xi\in TM$. Then for any $f\in C^{\infty} M$ we have a statistical structure $(g, D^f:=D+fK)$ on $M$ (see \cite{FHOSS}). In particular, if $M\simeq S^{2k+1}$ is a Sasakian sphere then for any $f\in C^{\infty} S^{2k+1}$ the collection ${(S^{2k+1}, g, D^f=D+fK)}$ is a statistical manifold of constant curvature 1.
\end{example}

\subsection{Statistical manifolds of constant curvature, dual connections, and affine immersions}\label{ku}

	
	\begin{defin}\label{da}
	Let $M\subset \R^{n+1}$ be an $n$-dimensional hypersurface and $\nabla$ be a standard connetion on $\R^{n+1}$. The section $\xi\ \in T\R^{n+1}|_M$ is a {\bfseries transversal vector field} along $M$ if 
	$$
	T\R^{n+1}|_M=TM\oplus\R \xi.
	$$
	When $\xi$ is given, we can define the {\bfseries induced affine connection} $D$ on $M$ and the {\bfseries second
	fundamental form} $h$ on $M$ as follows: for any $X,Y\in TM$ 
	$$
	\nabla_{X} Y =D_X Y+h(X,Y)\xi.
	$$
	If there exists a constant $\lambda$ such that $\xi=\lambda\left(\sum_{i=1}^{n+1} x^i\frac{\d}{\d x^i}\right)$ then the pair $\left(M,\xi\right)$ is called a {\bfseries centro-affine submanifold}.
	\end{defin} 

	\begin{defin}	
		Let $D$ be an affine connection on $M$, $\iota : M\to \R^n$ an immersion and $\xi \in \iota^*T\R^{n+1}$ such that $D$ is equal to the pullback of the induced affine connection $M$. Suppose for any neighborhood $U\subset M$ such that $\iota|_U$ is inclusion, the pair $(\iota U, \iota_* \xi)$ is a centro-affine submanifold. Then the pair $\left(\iota,\xi\right)$ is called a {\bfseries centro-affine immersion} of $(M,D)$.
	\end{defin}

	\begin{defin}
		Let $\left(M,g\right)$ be a Riemannian manifold. Two affine connection $D$ and $\bar D$ are called dual to each other (with respect to $g$) if 
		$$
		\L_X\left(g(Y,Z)\right)=g(D_X Y) + g(Y,\bar D_X Z).
		$$
	\end{defin}
	
	\begin{proposition}
		Let $D$ and $\bar D$ be dual affine connections on a Riemannian manifold $(M,g)$. Then $(M,g,D)$ is statistical if and only if both $\nabla$ and $\bar\nabla$ are torsion-free.
	\end{proposition}
	As a consequence, $(M,g,D)$ is a statistical manifold if and only if $(M,g,\bar D)$ is a statistical manifold.

	\begin{theorem}[\cite{K}]\label{K0}
		Let $\left(M,g,D\right)$ be a statistical manifold of dimension $d\ge 3$. If there exist centro-affine immersions $(\iota, \xi)$ of $(M,D)$ and $(\bar \iota, \bar \xi)$ of $(M,\bar D)$ such that the pullbacks of the second fundamental forms are equal to $g$ then $(M,g,D)$ is a statistical manifold of constant curvature $c$. Moreover, ${\nabla \xi = -c\text{Id}}$ and $\nabla \bar \xi= -c\text{Id}$, where $\nabla$ is the standard flat connection on $\R^n$.
	\end{theorem}
	\begin{theorem}[\cite{K}]\label{k}
		Let $\left(M,g,D\right)$ be a statistical manifold. If $\left(M,g,D\right)$ is a statistical manifold of constant curvature $c$ then there exist centro-affine immersions $(\iota, \xi)$ of $(M,D)$ and $(\bar \iota, \bar \xi)$ of $(M,\bar D)$ such that the second fundamental form is equal to $g$ and $\nabla \xi = -c\text{Id}$ and $\nabla \bar \xi= -c\text{Id}$, where $\nabla$ is the standard flat connection on $\R^n$. 
	\end{theorem}

	\subsection{Statistical manifolds of constant curvature and radiant Hessian manifolds}

	Consider a Hessian metric $g=\text{Hess} \ \varphi$ on $\R^n$ and the level set 
	$$
	M=\{x\in \R^{n+1} \ |\ \phi(x)=1\}.
	$$ 
	Suppose that for any $x\in M^{n+1}$ we have $d\phi_x\ne 0$. Then $M$ is a submanifold in $\R^{n+1}$. Let $E$ be the gradient of $\varphi$ with respect to the metric $g$. Then the field $E$ is transversal to $M$. 
	\begin{proposition}[\cite{Sh}, Lemma 5.1 and Example 5.1]\label{3.8}
		Let 
		$$\varphi,\ \ \ M=\{\varphi=1\}\subset\R^{n+1}, \ \ \ g=\text{Hess} \ \varphi,\ \ \ E=\text{grad} \ \varphi$$ be as above. Consider the induced affine connection $D$ on $M$ and the second fundamental form $h$ on $M$ (see Definition \ref{da}). Then 
		$$
		h=-\frac{1}{E(\phi)} g|_M.
		$$
		Moreover, suppose  $\left(M,E\right)$ is a centro-affine immersion and $\nabla E=\mu \text{Id}$. Then $\left(D,g|_M\right)$ is a statistical structure of constant curvature $c=\frac{\mu}{E(\varphi)}$.  
	\end{proposition}

			\begin{proposition}[\cite{Os}, Proposition 3.11]\label{311}
				Let $(C,g,\nabla,\xi)$ be a radiant Hessian manifold and $\nabla \xi=\lambda\text{Id}$, where $\lambda\in \R$. Then 
				$$
				g=\text{Hess}\left(\frac{g(\xi,\xi)}{4-2\lambda}\right).
				$$
			\end{proposition}

	\begin{theorem}\label{39}
		Let 	$
		\left(M\times \R^{>0}, g =s^2g_M+ds^2, \nabla, s\frac{\partial}{\partial s}\right) 
		$  
		be a radiant Hessian manifold and $\lambda\in \R$ a number which satisfies $\nabla\left( s\ds\right) =\lambda \text{Id}$. Then there exists a connection $D$ on $M$ such that $\left(M,g,D\right)$ is a statistical manifold of constant curvature $c=\lambda(2-\lambda)$.
	\end{theorem}
	
	\begin{proof}
		According to Proposition \ref{311}, the function $\varphi=\frac{s^2}{4-2\lambda}$ is a potential of $g$ and $M\times\{1\}$ is a level set of the potential. Then the gradient of $\phi$ is 
		$$
		E=\frac{1}{2-\lambda}s\ds.
		$$ 
		Since $\nabla \left(s \ds\right)=\lambda \text{Id}$, we have  
		$$
		\nabla E=\frac{\lambda}{2-\lambda}\text{Id}.
		$$
		Let $D$ be the affine connection induced by $E$. According to Proposition \ref{3.8}, $(M,g,D)$ is a statistical structure of curvature 
		$$
		c=\frac{\lambda}{(2-\lambda)E(\varphi)}=\lambda(2-\lambda).
		$$
	\end{proof}
\begin{rem} \label{r}
	Consider assumption  of Theorem \ref{39}. Let $h$ be the second fundamental form of the pair $(M,s\ds)$ from the proof of Theorem \ref{39}. By Proposition \ref{3.8},  we have
	$$
	h=-\frac{1}{E(\phi)} g=-{(2-\lambda)^2} g.
	$$
	Moreover, the vector field 
	$$
	\rho=\frac{1}{\lambda}s\ds
	$$
	on $M\times\ R^{>0}$ is radiant and
	$$
	E=\frac{\lambda}{2-\lambda}\rho.
	$$
	 
\end{rem}

	\begin{lemma}\label{ll}
		Let $\left(M\times\R^{>0},\nabla\right)$ be a flat affine manifold $s$ a coordinate on $\R^{>0}$, $\nabla \left(s\ds\right) =\lambda\text{Id}$, where $\lambda \in \R\setminus\{0,2\}$, and $\varphi= \left(\frac{s^2}{4-2\lambda}\right)$. Then $g=\text{Hess} \ \varphi$ can be written as
		$$
		g=s^2 g_M+ds^2,
		$$  
		where $g_M$ is a symmetric bilinear form on $M$.
	\end{lemma}
\begin{proof}
	We have
	\begin{multline*}
	\text{Hess}\ \varphi \left(s\ds, s\ds\right) =s\ds\left(s\ds\left(\frac{s^2}{4-2\lambda}\right)\right)-\left(\nabla_{s\ds} s\ds \right)\left(\frac{s^2}{4-2\lambda}\right)= \\
	= \frac{4s^2}{4-2\lambda}-\lambda s\ds \left(\frac{ s^2}{4-2\lambda}\right)=\frac{4s^2}{4-2\lambda}-\frac{2\lambda s^2}{4-2\lambda}=	
	s^2
	\end{multline*}
	
	and for any $X\in TM$
	$$
	\text{Hess} \varphi \left(s\ds, X\right)=s\ds\left(X(\varphi)\right)-\nabla_{s\ds} X(\varphi)=0
	$$
	because $X(\varphi)=0$ and $\nabla_{s\ds} X=\lambda X$. Moreover, 
	$$
	\L_{s\ds} \text{Hess}\left(\frac{s^2}{4-2\lambda}\right)=2\text{Hess}\left(\frac{s^2}{4-2\lambda}\right)
	$$ 
	Thus, 
	$
	g=s^2g_M+ds^2,
	$
	where $g_M$ is a symmetric bilinear form on $M$.
\end{proof}

	The following theorem is converse to Theorem \ref{39}.
	\begin{theorem}\label{38}
		Let $\left(M,g_M,D\right)$ be a statistical manifold of constant curvature $c\le1$ and $\lambda$ a solution of the equation $\lambda (2-\lambda)=c$ (as in Theorem \ref{39}). Then there exists a connection $\nabla$ on $M\times \R^{>0}$ such that
		$$
		\left(M\times \R^{>0}, g =s^2g_M+ds^2, \nabla, s\frac{\partial}{\partial s}\right) 
		$$  
		is a radiant Hessian manifold and ${\nabla \left(s\ds\right)=\lambda \Id}$.
	\end{theorem}
	
	\begin{proof}
		We proved this theorem in \cite{Os3} by direct calculations. Here, we present the proof based on centro-affine immersions. Let $\lambda$ a solution of the equation $\lambda (2-\lambda)=c$ (as in Theorem \ref{39}). Then $(-\left(2-\lambda)^2 g,D\right)$ is a statistical structure of curvature $-{\frac{\lambda}{2-\lambda}}$.  
		According to Theorem \ref{k}, there exists a centro-affine immersion ${\iota :M \to \R^{n+1}}$ with a transversal field $\xi=\frac{\lambda}{2-\lambda}\sum_{i=1}^{n+1}x^i\frac{\d}{\d x^i}$ and the second fundamental form $-{(2-\lambda)^2} g$ (as in Remark \ref{r}). 
		Set an immersion 
		$$
		\hat \iota: M\times \R^{>0} \to \R^{n+1}, \ \ \ \ \hat \iota (m\times s)=s^{\lambda} \iota(m).
		$$
		Then we have $s\ds=\hat\iota^* \left(\lambda\sum_{i=1}^{n+1}x^i\frac{\d}{\d x^i}\right)$.  Set a connection $\nabla$ on $M\times \R^{>0}$ as the pullback of the standard connection on $\R^{n+1}$ and 
		$$
		\widetilde g= \text{Hess}\frac{s^2}{4-2\lambda}.
		$$ 
		According to Lemma \ref{ll}, $\widetilde g$ can be written as 
		$$
		\widetilde g=s^2 \widetilde g_M+ds^2.
		$$
		The vector 
		$E=\frac{1}{2-\lambda} s\ds$
		is the gradient of $\frac{s^2}{4-2\lambda}$ (as in Remark \ref{r}).
		According to Proposition \ref{3.8}, 
		$$
		\widetilde g_M=-E(\varphi) h =g_M.
		$$
		Thus, $\left(M\times\R^{>0},s^2g_M+ds^2,\nabla,s\ds\right)$ is a radiant Hessian manifold.
	\end{proof}

	We constructed a correspondence between radiant Hessian manifolds and statistical manifolds of constant curvature.	

	\section{Locally conformally Hessian manifolds}
	\subsection{Locally conformally Hessian manifolds: definition and examples}

	\begin{defin}
		A {\bfseries locally conformally Hessian (l.c.H.) manifold} $(C,\nabla, g,\theta)$ is a flat affine manifold $(C,\nabla)$ endowed with a Riemannian metric $g$ and a closed 1-form $\theta$ such that $\nabla g-\theta\otimes g$ is a totally symmetric tensor. The form $\theta$ and the vector field $\xi =	\theta^\sharp$ are called {\bfseries Lee form} and {\bfseries Lee vector field}. 
	\end{defin}

	Suppose that $(C,\nabla, g,\theta)$ is an l.c.H. manifold. The closed form $\theta$ is locally defined as a differential of a function $\theta=df$. Since we have
	$$
	\nabla \left(e^{-f} g\right) = e^{-f}\nabla g -e^{-f}df\otimes g,
	$$
	the tensor $\nabla \left(e^{-f} g\right)$ is totally symmetric. Hence, the locally defined metric $e^{-f}g$ is Hessian (see Definition \ref{Hes}).
	
	We are interested in l.c.H. manifold with an affine Killing Lee vector field.

	\begin{example}
		A {\bfseries Hopf manifold} is $H^n=\R^n\setminus\{0\} /_	{x\sim a Ax}$, where $a\in \R^{>0}\setminus \{1\}$ and $A\in \text{O}(n)$. The universal covering $\R^n\setminus \{0\}$ has the l.c.H. metric $\frac{\sum \left(dx^i\right)^2}{\sum \left(x^i\right)^2}$ which is invariant with respect to homotheties. Hence, this metric induces an l.c.H. structure on the Hopf manifold. The Lee vector field is $\sum x^i \frac{\d}{\d x^i}$. 
	\end{example}
	\begin{example}
		Poincaré half-space $\left(\R^{>0}\times\R^n,\frac{\sum_{i=0}^n \left(dx^i\right)^2}{\left(x^0\right)^2}\right)$ is an l.c.H. manifold. The Lee vector field equal
		$x^0\frac{\partial}{\partial x^0}$. This field is affine but not Killing. Consider the action of $\Z^{n+1}$ on $\R^{>0}\times \R^n$ defined by 
		$$
		\left(z^0,z^1,\ldots,z^n\right) \left(x^0,x^1,\ldots, x^n\right)=  \left(\lambda^{z^0}x^0, x^1+z^1,\ldots,x^n+z^n\right),
		$$ 
		where $\lambda\in \R^{>0}\setminus \{1\}$. This action preserves the affine structure and the metric $\frac{\sum_{i=1}^n \left(dx^i\right)^2}{\left(x^0\right)^2}$. Hence, this metric induces an l.c.H. structure on the torus $\R^{>0}\times\R^n/\Z^n =T^{n+1}$.
	\end{example}

\begin{example}
	Let $V$ be a convex cone without any straight full line, $\text{Aut}(V)$ the group of linear automorphisms of $V$, and $\psi$ the characteristic function (see Example \ref{e33}). Then 
	$$
	g_{H.}=\text{Hess} \ln \psi \ \ \ \ \text{and} \ \ \ \ g_{l.c.H.}=\frac{\text{Hess}\ \psi}{\psi}
	$$	
	are $\text{Aut}(V)$-invariant Hessian and l.c.H. forms (\cite{V}). The Lee form is equal to $-d\ln \varphi=-\frac{d\varphi}{\varphi}$. Then $\frac{\text{Hess}\ \psi}{\psi}$ induces an l.c.H. structure on $V/\Gamma$. Suppose $\Gamma$ is a subgroup of linear automorphisms $\text{Aut}(V)$ such that $V/\Gamma$ is a compact manifold. Then $g_{H.}$ and $g_{l.c.H.}$ induces Hessian and l.c.H. metrics on $V/\Gamma$. Note that this situation differs from the K\"ahler case: a compact manifold cannot admits K\"ahler and  l.c.K. structure at the same time (\cite{Va}). 
\end{example}
	\begin{example}\label{e67}
	Consider an affine plane $\R^2$ with affine coordinates $(x,y)$ and conformally Hessian metric 
	$$
	g=\frac{\text{Hess} \left(e^x+e^y\right)}{(1+e^{ \frac{x}{2} })^2+(1+e^{\frac{y}{2}} )^2}=\frac{e^xdx^2+e^ydy^2}{(1+e^{ \frac{x}{2} })^2+(1+e^{\frac{y}{2}} )^2}.
	$$
	Let us check that the Lee vector field is Killing but not affine. Set $\tilde x=1+e^{\frac{x}{2}}$ and $\tilde y=1+e^{\frac{y}{2}}$. Then 
	$$
	g=\frac{4 d\tilde x ^2+4 d\tilde y ^2}{\tilde x^2+\tilde y^2}.
	$$
	The Lee form is
	$$
	d\ln\left(\frac{1}{\tilde x^2+\tilde y^2} \right)=
	\frac{
		-2\tilde xd\tilde x-2\tilde y d \tilde y
	}
	{\tilde x^2+\tilde y^2}.
	$$
	The Lee vector field is 
	$$
	\xi=-\frac{\tilde x}{2}\frac{\d}{\d \tilde x}-\frac{\tilde y}{2}\frac{\d}{\d \tilde y}.
	$$
	This vector field is Killing. In flat affine coordinates we have 
	$$
	\xi=-\frac{1+e^{\frac{x}{2}}}{2}\frac{\partial}{\partial \left(1+e^{\frac{x}{2}}\right)}-\frac{1+e^{\frac{y}{2}}}{2}\frac{\partial}{\partial \left(1+e^{\frac{y}{2}}\right)}=-\frac{\left(e^{\frac{x}{2}}+e^x\right)\frac{\d}{\d x}+\left(e^{\frac{y}{2}}+e^y\right)\frac{\d}{\d y}}{(1+e^{ \frac{x}{2} })^2+(1+e^{\frac{y}{2}} )^2}.
	$$ 
	Therefore the Lee vector field $\xi$ is Killing but not affine.
\end{example}

\begin{proposition}
	Let $(C,\nabla,g,\theta)$ be a compact l.c.H manifold and $\theta\ne 0$ at any point. Then $(g,\nabla)$ is not globally conformally Hessian.
\end{proposition}	

\begin{proof}
	Suppose that $g$ is globally conformally Hessian. That is, there exists a function $p$ on $C$ such that the tensor
	$$
	\nabla (p g) =p \nabla g+dp\otimes \theta
	$$
	is totally symmetric. Since $C$ is compact, there exists a point $x\in C$ such that $dp|_{T_x C}=0$. Hence, 
	$$\nabla g|_{T_x C}=0
	$$ is totally symmetric. Since $\left(C,\nabla,\theta\right)$ is an l.c.H. manifold, the tensor
	$
	\nabla g-\theta\otimes g.
	$
	is totally symmetric. Thus, 
	$$
	\theta\otimes g|_{T_x C}=0
	$$ 
	is totally symmetric. Choose $X,Y\in {T_x C}$, such that $\theta(X)=0$ and $\theta(Y)\ne 0$. Then 
	$$
	0\ne \theta(Y)\otimes g(X,X)=\theta(X)\otimes g(X,Y)=0.
	$$ 
	We get the contradiction. Thus, $g$ is not globally conformally Hessian.
\end{proof}

The examples of statistical manifolds of constant curvature from {Section \ref{s4}} combining with the following theorem provide examples of l.c.H. manifolds.
		\begin{theorem}\label{65}
		Let $\phi$ be an automorphism of a statistical manifold $(M,g,D)$  of constant curvature $c\le 1
		$, $\lambda$ a solution of the equation $c=\lambda(2-\lambda)$, and $q\in \R^{>0}$. Consider the automorphism 
		$$
		{\phi_q: M\times \R^{>0}\to M\times \R^{>0}, \ \ \ \phi_q(m,t)= \left(\phi(m),qt\right)}.
		$$
		Then there is a connection $\nabla$ on $M\times \R^{>0}/_{\phi_q}$ such that 
		$$
		\left(M\times \R^{>0}/_{\phi_q}, \ g_M+\frac{ds^2}{s^2}, \ \frac{-2ds}{s}\right)
		$$
		is an l.c.H. manifold with an affine Killing Lee vector field $\xi=-2s\ds$ satisfying ${\nabla \xi =-2\lambda \Id.}$
	\end{theorem}
\begin{rem}
	The vector field $s\ds$ on $M\times \R^{>0}$ and the differential form $\frac{ds}{s}$ on $M\times \R^{>0}$ are invariant with are invariant with respect to $\varphi_q$. Hence, we can consider $s\ds$, $\frac{ds}{s}$, and $\frac{ds^2}{s^2}$ as tensors on $M\times \R^{>0}/_{\phi_q}$.
\end{rem}

	\begin{proof}[Proof of Theorem \ref{65}.]
			According to Theorem \ref{38}, there exists a connection $\nabla$ on $M\times \R^{>0}$ such that the collection 
		$
		\left(M\times \R^{>0}, s^2g_M+ds^2, \nabla, s\frac{\partial}{\partial s}\right) 
		$  
		is a radiant Hessian manifold. The automorphism $\phi_q$ preserves the metric $g_M+\frac{ds^2}{s^2}$. Moreover, since the field $s\ds$ is affine, $\phi_q$ preserves the affine structure $\nabla$. Therefore, the pair $\left(\nabla,g_M+\frac{ds^2}{s^2}\right)$ induces an l.c.H. structure on $M\times \R^{>0}/_{\phi_q}$. The Lee form equals
		$$
		d\ln \frac{1}{s^2}=\frac{-2ds}{s}
		$$
		and the Lee vector equals $-2s\ds$. This field is affine and Killing.
		
		Therefore, 
		$
		\left(M\times \R^{>0}/_{\phi_q}, \ g_M+\frac{ds^2}{s^2}, \ \frac{-2ds}{s}\right)
		$
		is a l.k.H. manifold with an affine Killing Lee vector field.
	\end{proof}
	\begin{proof}[Proof of Theorem \ref{t1}]
		If $(M,g_M,D)$ is a statistical manifold of curvature $c>1$, we replace $g_M$ by $cg_M$ and get a statistical manifold of curvature 1. Using Theorem \ref{65}, we can construct a l.c.H. structure with an affine Killing Lee vector on $M\times \R^{>0}/_{\phi_q}$.
	\end{proof}

	\subsection{Minimal Hessian covering}
	\begin{defin}
		Let $\left(C,g,\nabla, \theta\right)$ be an l.c.H. manifold. Then the {\bfseries weight bundle} is the trivial $\R$-bundle $L$ endowed with a flat connection defined by  $\theta$.
	\end{defin}
	\begin{defin} \label{3.3}
	Let $(C,\nabla, g,\theta)$ be an l.c.H. manifold. Fix a point $x_0\in C$. Consider the functional $L:\pi_1(C,x_0) \to \R$ defined by
	$$
	 L(\gamma)=\int_\gamma \theta.
	$$
	Then the image of $L$ is the monodromy group $\Gamma$ of the weight bundle $(L,\theta)$.
	Let $\widetilde{C}$ be the covering of $C$ with the fiber $\Gamma$. For any $\gamma\in \pi_1(\widetilde C)$ we have 
	$$
	\int_\gamma \pi^*\theta=0.
	$$
	Therefore, there exists a function $f$ on $\widetilde C$ such that $df=\pi^*\theta$. The metric $\widetilde g=e^{-f}\pi^*g$ is Hessian. We say that that $\left(\widetilde{C},\pi^* \nabla,\widetilde g\right)$ is the {\bfseries minimal Hessian covering} of the l.c.H. manifold $(C,\nabla, g,\theta)$. The deck transform group of the covering $ \widetilde C \to C$ act on $\widetilde {C}$ by Hessian homotheties. 	
	\end{defin}
	\begin{proposition}[\cite{OV4}]\label{3.4}
		Let $\theta$ be a 1-form on a Riemannian manifold $(C, g)$, ${\xi=\theta^\sharp}$, and $\nabla_\text{LC}$ the Levi-Civita connection. Then the following conditions are equivalent:
		\begin{itemize}
			\item [(i)] $\nabla_\text{LC} \xi=0$.
			\item [(ii)]  $\nabla_\text{LC}\theta=0$.
			\item [(iii)] $d\theta=0$ and $\L_\xi g = 0$.
		\end{itemize}
	\end{proposition}
	\begin{cor}\label{69}
		Let $(C,\nabla, g,\theta)$ be an l.c.H. manifold with a Killing Lee vector $\xi$. Then $a:=\theta(\xi)=g(\xi,\xi)$ is constant.  
	\end{cor}
\begin{proof}
	According to Proposition \ref{3.4}, $\nabla_\text{L-C} \xi=0$. Therefore, $g(\xi,\xi)$ is constant.
\end{proof}

\begin{lemma}\label{61}
	Let $(C,\nabla, g,\theta)$ be an l.c.H. manifold, $\xi$ a vector field on $C$, $\left(\widetilde C,\pi^* \nabla,\widetilde g\right)$ the minimal Hessian covering (where $\pi:\widetilde{C}\to C$ is the corresponding covering map), and $a=g(\xi,\xi)=\theta(\xi)$. Then $\L_\xi g=0$ if and only if $\left( \widetilde C,\widetilde g,\widetilde{\xi}:=-\frac{2}{a}\pi^*\xi\right)$ is
	a self-similar manifold. 
\end{lemma}

\begin{proof}
	There exists a function $f$ on $\widetilde{C}$ such that $\pi^*\theta =df$ and $\widetilde g=e^{-f}\pi^* g$ (see Definition \ref{3.3}). Therefore, 
	$$
	\widetilde{\xi}(f)=\left\langle\widetilde{\xi},\pi^*\theta \right\rangle =\left\langle-\frac{2}{\theta(\xi)}\pi^*\xi,\pi^*\theta \right\rangle= -2
	$$
	Hence,
	$$
	\L_{\widetilde{\xi}} \widetilde g= \L_{\widetilde{\xi}}\left(e^{-f}\pi^* g\right)=-e^{-f}\widetilde\xi(f)\pi^*g+e^{-f}\L_{\widetilde{\xi}} \pi^*g=2\widetilde g+e^{-f}\pi^*\left(\L_\xi g\right).
	$$
	Thus, $\L_{\widetilde{\xi}} \widetilde{g}=2\widetilde{g}$ if and only if $\L_\xi g =0$.
\end{proof}
 
\begin{proposition}\label{35}
	Let $(C,\nabla, g,\theta)$ is an l.c.H. manifold with a Killing Lee vector $\xi$, $a=g(\xi,\xi)$, and  $\left(\widetilde C,\pi^* \nabla,\widetilde g\right)$ its minimal Hessian covering.
	Then $\left( \widetilde C,\widetilde g,\widetilde{\xi}:=-\frac{2}{a}\pi^*\xi\right)$ is a self-similar manifold with a potential homothetic vector field. 
\end{proposition}

\begin{proof}
	According to Lemma \ref{61}, $\left( \widetilde C,\widetilde g,-\frac{2}{a}\pi^*\xi\right)$ is a self-similar manifold.
	
	There exists a function $f$ on $\widetilde{C}$ such that $\pi^*\theta =df$ and $\widetilde g=e^{-f}\pi^* g$ (see Definition \ref{3.3}). 	Since $\iota_\xi g= \theta$, we have 
	$$
	\iota_{\pi^*\xi} \pi^* g =\pi^* \theta =df
	$$
	and
	$$
	\iota_{\pi^*\xi} \widetilde g =e^{-f} df= d\left(-e^{-f}\right).
	$$
	Thus the field ${\pi^*\xi}$ is potential. According to Corollary \ref{69}, $g(\xi,\xi)$ is constant. Therefore, the homothetic field $-\frac{2}{a}\pi^*\xi$ is potential.
\end{proof}
\subsection{Radiant l.c.H. manifolds}
\begin{defin}\label{414}
	Let $(C,\nabla, g,\theta)$ be an l.c.H. manifold with a Killing Lee vector field $\xi$, $a=g(\xi,\xi)$, and  $\left(\widetilde C,\pi^* \nabla,\widetilde g\right)$ its minimal Hessian covering.
	We say that $(C,\nabla, g,\theta)$ is a {\bfseries radiant l.c.H} manifold if  $\left( \widetilde C,\widetilde g,\widetilde \xi=-\frac{2}{a}\pi^*\xi\right)$ is a radiant Hessian manifold. 
\end{defin}

\begin{proposition}\label{415}
	Let  $(C,\nabla, g,\theta)$ be a l.c.H. manifold with a Killing Lee vector field and $a=g(\xi,\xi)$. Then $(C,\nabla, g,\theta)$ is a radiant l.c.H. manifold if and only if there exists a constant $\mu\in\R\setminus\{-a,0\}$ such that ${\nabla \xi=\mu\Id}$.
\end{proposition}
\begin{proof}
	Suppose there exist a constant $\mu$ such that ${\nabla \xi=\mu\Id}$. According to Proposition \ref{35}, $\left( \widetilde C,\widetilde g,\widetilde \xi=-\frac{2}{a}\pi^*\xi\right)$ is a self-similar manifold. Then 
	$$
	\pi^*\nabla{\widetilde \xi}= \lambda\Id,
	$$
	 where 
	 $$
	 \lambda=-\frac{2}{a}\mu
	 $$ 
	 is a constant. The proof in the opposite direction is similar. We have $\lambda \in \R\setminus \{0,2\}$ (see Remark \ref{rem}). Therefore, $\mu\in \R \setminus \{-a,0\}$.
\end{proof}

\begin{theorem}\cite{Ga} \label{Ga}
	Let a group $\Gamma$ act on a Riemannian cone $$\left(M\times\R^{>0},s^2g_M+ds^2\right)$$ 
	by homotheties and let the quotient $M\times\R^{>0}/\Gamma$ be a compact manifold. Then the holonomy of
	 $\left(M\times\R^{>0},s^2g_M+ds^2\right)$ is irreducible or $\left(M\times\R^{>0},s^2g_M+ds^2\right)$ is flat.
\end{theorem}

\begin{proof}[Proof of Theorem \ref{TT3}]
	According to Proposition \ref{35}, $\left(\widetilde C,\widetilde g,\widetilde \xi\right)$ is a self-similar manifold with potential homothetic vector field $\widetilde \xi$. Since the manifold ${C}$ is compact the vector field $\xi$ is complete. Hence, the vector field $\widetilde \xi$ is complete. Combining it with Theorem \ref{T1}, we get that $\left(\widetilde C,\widetilde g\right)$ is isometric to the Euclidean space with a radiant vector field $\left(\R^n,\sum_{i=1}^n\left(dx^i\right)^2,\sum x^i\frac{\d}{\d x^i} \right)$ or a to Riemannian cone.
	 
	Suppose $\left(\widetilde{C},\widetilde{g},\widetilde{\xi}\right) =\left(\R^n,\sum_{i=1}^n\left(dx^i\right)^2,\sum x^i\frac{\d}{\d x^i} \right)$. Let $\Gamma$ be the monodromy group of the covering $\pi: \widetilde{C}\to C$. Then $\Gamma$ acts on $\left(\R^n,\sum_{i=1}^n\left(dx^i\right)^2 \right)$ by homotheties preserving $\sum x^i\frac{\d}{\d x^i}$. Choose an element $\gamma\in \Gamma$ such that the action is not an isometry. Then the action of $\gamma$ or the action of $\gamma^{-1}$ is a contraction map preserving the point $0\in\R^n$. Hence, the quotient $\widetilde{C}/\Gamma$ cannot be Hausdorff. However, $\widetilde{C}/\Gamma=C$ is a manifold.
	
		Thus, $\left(\widetilde C,\widetilde g\right)$ is isometric to a Riemannian cone. According to Theorem \ref{Ga}, the holonomy of $\left(\widetilde C,\widetilde g\right)$ is irreducible or $\left(\widetilde C,\widetilde g \right)$ is flat. Since $(C,g)$ is not locally conformally flat, $\left(\widetilde C,\widetilde g \right)$ cannot be flat. Hence, the holonomy of $\left(\widetilde C,\widetilde g\right)$ is irreducible. According to Corollary \ref{c1}, $\left(\widetilde C,\pi^*\nabla,\widetilde g,\widetilde\xi\right)$ is a radiant Hessian manifold. 
\end{proof}

	\subsection{An l.c.H. metric expressed in terms of the Lee form}\label{s8}
	\begin{proposition}\label{71}
		Let $\left(C,\nabla, g,\theta\right)$ be a compact radiant l.c.H. manifold $\xi$, the Lee vector field, $\nabla \xi=\mu \Id$, and $a=g(\xi,\xi)$ the constant from Corollary \ref{69}. Set $u=-\mu-a$. Then we have 
		$$
		ug=\nabla \theta  - \theta\otimes \theta,
		$$
		where $u\ne0$.
	\end{proposition}

\begin{proof}
	 We have 
	$$
	\left(\nabla g \right)(\xi,X,Y) = \L_\xi \left(g(X,Y)\right)- g(\nabla_\xi X, Y)- g(X,\nabla_\xi  Y).
	$$
	Since $\xi$ is Killing,
	$$
	\L_\xi (g(X,Y)) =g([\xi, X], Y)+ g(X,[\xi, Y]).
	$$
	Moreover, the connection $\nabla$ is torsion free. Thus,
	$$
	\left(\nabla g \right)(\xi,X,Y)=-g(\nabla_X \xi,Y)-g(X,\nabla_Y \xi).
	$$
	According to Theorem \ref{TT3}, there exists	 $\mu \in \R$  such that $\mu\notin \{0,-a\}$ and $\nabla \xi=\mu\Id$. Therefore, 
    $$
	\left(\nabla g \right)(\xi,X,Y)=-g(\nabla_X \xi,Y)-g(X,\nabla_Y \xi)=-{2\mu } g(X,Y).
	$$
	and 
	\begin{equation}\label{eee1}
	\left(\nabla g -\theta \otimes g\right)(\xi,X,Y)=(-2\mu - a) g(X,Y).
	\end{equation}
	Using the identities $\iota_\xi g =\theta$ and $\nabla_X\xi=\mu \text{Id}$ we obtain
	\begin{multline}\label{eee2}
	(\nabla g -\theta\otimes g)(X,\xi,Y)=\L_X(\theta(Y))-\theta(\nabla_X Y)-g(\nabla_X \xi, Y)-\theta (X)\theta(Y)= \\
	{-\mu}{g(X,Y)}+(\nabla\theta -\theta\otimes \theta) (X,Y).
	\end{multline}
	By the definition of l.c.H. manifolds, the tensor $\nabla g-\theta\otimes g$ is totally symmetric. Combining this with \eqref{eee1} and \eqref{eee2}, we get that
	$$
	\left(-\mu-a\right)g=\nabla\theta -\theta\otimes \theta.
	$$
	According to Corollary \ref{69}, we have $\mu\ne -a$. Hence, $u=-\mu-a\ne 0$.
\end{proof}

\begin{proposition}\label{7.2}
	Let $(C,\nabla)$ be a flat affine manifold and $\theta$ a closed 1-form such that the bilinear form 
	$
	\nabla \theta -\theta\otimes \theta 
	$
	is positive definite. Then $(C,\nabla,g_\theta=\nabla \theta -\theta\otimes \theta,\theta)$ is an l.c.H. manifold.  
\end{proposition}

\begin{proof}
	It is enough to check that $(C,\nabla,g_\theta,\theta)$ is an l.c.H. manifold locally. The form $\theta$ is locally expressed as $\theta = df$. Then 
	$$
	\nabla \left(d\left(-e^{-f}\right)\right)=\nabla \left(e^{-f} df\right)=e^{-f}\nabla df -e^{-f} df\otimes df=e^{-f} g_\theta. 
	$$
	Thus the metric $e^{-f} g_\theta$ is Hessian and $g$ is l.c.H.
\end{proof}

\begin{proposition}\label{p00}
	Let $u$ be a nonzero constant, $$(C,\nabla, g=u^{-1}(\nabla\theta-\theta\otimes \theta),\theta)$$ a  l.c.H. manifold and $\xi$ be a Killing vector field on $(C,g)$ such that $\L_\xi \theta=0$ and $\nabla \xi=\mu \Id$, for a constant $\mu\in \R$. Then $\xi$ coincides with the Lee vector field up to a constant  multiplier.
\end{proposition}

\begin{proof}
	Since 
	$
	\L_\xi \theta=0,
	$
	and $\theta$ is closed, the value $\theta(\xi)$ is constant. Hence, for any $X\in TC$
	$$
	\left(\nabla \theta\right)(\xi, X)=\L_X\left(\theta(\xi)\right)-\theta\left(\nabla_X \xi\right)=-\mu \theta(X).
	$$
	Therefore,
	$$
	\iota_\xi g=\frac{\iota_\xi (\nabla \theta -\theta\otimes \theta)}{u}=\frac{(-\mu-\theta(\xi))\theta}{u},
	$$
	i.e. $\iota_\xi g$ is proportional to $\theta$. Thus, $\xi$ coincides with the Lee vector field up to multiplication on a constant.
\end{proof}

Combining Propositions \ref{71} and \ref{p00} we get the following.

\begin{theorem}\label{t85}
	Let $(C,\nabla, g,\theta)$ a compact radiant l.c.H. manifold and $\xi$ be a Killing vector field on $(C,g)$ such that $\L_\xi \theta=0$ and $\nabla \xi=\mu \Id$, for a constant $\mu\in \R$. Then $\xi$ coincides with the Lee vector field up to a constant  multiplier.
\end{theorem}

\begin{defin}[\cite{Sh}]
	A Hessian manifold $(C,\nabla, g)$ is said to be of {\bfseries Koszul type} if there exists a (globally defined) closed 1-form $\theta$ such that $g=\nabla\theta$. 
\end{defin}

\begin{theorem}[\cite{Ko}]\label{K}
	Let $(C,\nabla, g)$ be a compact Hessian manifold of Koszul type. Then the universal covering $\widetilde{C}$ is a convex cone without full straight lines. Moreover, the lifting of the Hessian metric equals to $\text{Hess} \left(\ln \psi\right)$ up to a constant  multiplier, where $\psi$ is the characteristic function of the cone (see example \ref{e33}). 
\end{theorem}

\begin{theorem}\label{t6}
		Let $\left(C,\nabla, g,\theta\right)$ be a compact radiant l.c.H. manifold, $\xi$ the Lee vector field, $\nabla \xi =\mu \Id$, and $a=g(\xi,\xi)$. Suppose $-\mu-a>0$. Then $\left(C,\nabla,\nabla\theta\right)$ is a Hessian manifold of Koszul type.
\end{theorem}
\begin{proof}
	According to Proposition \ref{71}, we have 
	$$
	ug=\nabla\theta -\theta\otimes\theta,
	$$
	where $u=-\mu-a>0$. Therefore, $\nabla \theta$ is positive definite and $\left(C,\nabla,\nabla \theta\right)$ is compact Hessian manifold of Koszul type.
\end{proof}

Combining Theorems \ref{t6} and \ref{K} we get the following.

\begin{cor}
	Let $\left(C,\nabla, g,\theta\right)$ be a compact radiant l.c.H. manifold, $\xi$ the Lee vector field, $\nabla \xi =\mu \Id$, and $a=g(\xi,\xi)$. Suppose $-\mu-a>0$.  Then
	the universal covering of $C$ is a convex cone without full straight lines.
\end{cor}

\begin{theorem}\label{K1}
Let $\phi$ be an automorphism of a statistical manifold $(M,g,D)$  of constant curvature $c<0$. Consider the automorphism 
$$
{\phi_q: M\times \R^{>0}\to M\times \R^{>0}, \ \ \ \phi_q(m,t)= \left(\phi(m),qt\right)}.
$$
Then $M\times \R^{>0} /_{\phi_q}$ admits a Hessian structure of Koszul type.
\end{theorem}

\begin{proof}
	Let $\lambda$ be the positive solution of the equation $c=\lambda(2-\lambda)$. $$
	{\phi_q: M\times \R^{>0}\to M\times \R^{>0}, \ \ \ \phi_q(m,t)= \left(\phi(m),qt\right)}.
	$$
	According to Theorem \ref{65}, there is a connection $\nabla$ on $M\times \R^{>0}/_{\phi_q}$ such that 
	$$
	\left(M\times \R^{>0}/_{\phi_q}, \nabla,\ g=g_M+\frac{ds^2}{s^2}, \ \theta=\frac{-2ds}{s}\right)
	$$
	is a radiant l.c.H. manifold with Lee vector field $\xi=-2s\ds$ satisfying ${\nabla \xi = -2\lambda}$. 

	Since $\lambda>2$, the constant $\mu=-2\lambda$ and $a=g(\xi,\xi)=4$ satisfies the condition of Theorem \ref{t6}. Therefore, $\left(M\times \R^{>0} /_{\phi_q},\nabla, \nabla \theta\right)$ is Hessian manifold of Koszul type. 
\end{proof}

\begin{cor}\label{C1}
	Let $\left(M,g_M,D\right)$ be a compact statistical manifold of negative constant curvature. Then the universal covering $\left(\widetilde M,\widetilde {g_M}\right)$ is isometric to the characteristic hypersurface of a cone.
\end{cor}
\begin{proof}
	According to Theorems \ref{K1}, there is a Hessian structure of Koszul type
	 $\left(\nabla,g=\nabla\theta\right)$ on 
	$M\times \R^{>0} /_{\phi_q}$, where 
	$\theta=\frac{-2ds}{s}$. According to Theorem \cite{K}, the universal covering $\left(\widetilde{M}\times \R^{>0}, \widetilde{g}\right)$ is isometric to convex cone without straight lines with the metric $\text{Hess}\left(\ln \psi\right)$, where $\psi$ is the characteristic function (see Example \ref{e33}). Hence, $\theta=d \ln \psi$. The set $\widetilde M\times 1$ is an integral hypersurface of the 1-form $\theta$. That is, $\widetilde M\times 1$ is a level set of the function $\psi$. Therefore, $\left(\widetilde{M},\widetilde {g_M}\right)=\left(\widetilde{M},\widetilde {g}|_{\widetilde {M}}\right)$ is the characteristic hypersurface of the cone $\widetilde{M}\times \R^{>0}$.
\end{proof}

We can identify a characteristic hypersurface of a cone with a properly convex domain in $\mathbb{RP}^n$ (see Example \ref{e34}). Therefore, Corollary \ref{C1}  is equivalent to the second part of Theorem \ref{OT}.

\subsection{L.c.H manifolds of rank 1}

\begin{defin}
	We say that an l.c.H. manifold $\left(C,\nabla,g,\theta\right)$ is of of rank 1 if the monodromy group of the weight bundle $(L,\theta)$ is isomorphic to $\Z$.
\end{defin}

\begin{proposition}\label{426}
	 Let $\left(C,\nabla,g,\theta\right)$ be a l.c.H. manifold. The monodromy group of the weight bundle $(L,\theta)$ is isomorphic to $\Z$ if and only if ${[\theta]\in H^1(C,\Q)}$.
\end{proposition}
\begin{proof}
	The proof coincides with the proof of the similar proposition for l.c.K. manifolds in \cite{OV3}.
\end{proof}
\begin{proof}[Proof of Theorem \ref{14}]
	Consider a cohomology class $[\alpha] \in H^1(M;\R)$ and let $\alpha$ be its harmonic representative. According to Proposition \ref{71}, $g$ is proportional to $\nabla\theta-\theta\otimes\theta$. If $\alpha$ is chosen sufficiently small than $g_{\theta'}=\nabla \theta'-\theta'\otimes\theta'$ is positive definite. According to Proposition \ref{7.2}, $\left(C,\nabla,g_{\theta'},\theta'\right)$ is an l.c.H. manifold.
	
	Let us show that the Lee vector field $\xi$ of $(C,\nabla,g,\theta)$ coincides with the Lee vector field of $\left(C,\nabla,g_{\theta'},\theta'\right)$ up to a constant multiplier. Since the flow along the Lee vector field $\xi$ acts on $(C,g)$ by isometries and the form $\alpha$ is harmonic, this flow preserves $\alpha$. Hence it preserves the form $\theta'$. Moreover the flow along $\xi$ preserves $\nabla$, the flow along $\xi$ preserves $g_{\theta'}=\nabla \theta'-\theta'\otimes\theta'$. According to Proposition \ref{p00}, there exists a constant $a\in \R$ such that $a\xi$ is the Lee vector field on  $(C,\nabla,g',\theta')$. Since the field $\xi$ is affine and Killing with respect to $g_{\theta'}$, the field $a\xi$ is affine and Killing too. Thus, $(C,\nabla,g',\theta')$ is a radiant l.c.H. manifold of l.c.H. rank 1.
\end{proof}

Combining Theorem \ref{14} and Proposition \ref{426}, we get the following.

\begin{theorem}
	Let $(C,\nabla,g,\theta)$ be a compact radiant l.c.H. manifold. Then $g$ can be approximated by a sequence of Riemannian metrics which are conformally
	equivalent to an l.c.H. metric of rank 1 with an affine Killing Lee vector.
\end{theorem}


\begin{proof}[Proof of Theorem \ref{t2}]
	The minimal Hessian covering $(\widetilde C, \widetilde \nabla,\widetilde g)$ endowed with the vector field $\widetilde{\xi}=-\frac{2}{a}\pi^* \xi$ is a radiant Hessian manifold. We have 
	$$
	\widetilde{\nabla} \widetilde \xi= \lambda, \ \ \ \text{where} \  \ \ \lambda=-\frac{2}{a}.
	$$
	By Corollary \ref{69}, $a=g(\xi,\xi)$ is constant. Hence, $\xi$ is nonzero at any point. According to Theorem \ref{39}, $(\widetilde C, \widetilde g)$ is isometric to a cone $\left(M\times\R^{>0},s^2 g_M+ds^2\right)$ over a statistical manifold $(M,g_M,D)$ of constant curvature $c=\lambda(2-\lambda)$. Since $\lambda=-\frac{2}{a}$, the constant $c$ is negative if and only if $\mu\in (-\infty, -a)\cup (0,\infty)$.
	
	Since $M$ is an l.c.H. manifold of rank $1$, the deck group $\Gamma$ of $\widetilde C$ is isomorphic to $\Z$. The manifold $C$ is obtained from $\widetilde C$ as a factor $C=\widetilde{C}/\Gamma$, where $\Gamma$ acts on $\widetilde C =M\times \R^{>0}$ be homotheties. Let the generator $\gamma$ of $\Gamma$ act on the first component of $M\times \R^{>0}$ by an isometry $\varphi$ and on the second component  by multiplication by $q\in \R^{>0}$. Set a submersion 
	$$
	\sigma_0\ : \ C\simeq (M\times \R^{>0})/\Gamma\ \to \ \R^{>0}/\{a\}\simeq S^1
	$$
	using the diagram
		\begin{equation*}
		\begin{CD}
		M\times\R^{>0}@>\sigma>>  \R^{>0}  \\
		@V\pi VV  @V\pi_0VV              \\
		 (M\times \R^{>0})/\Gamma @>{\sigma_0}>>  \R^{>0}/\{q\} 
		\end{CD}.
		\end{equation*}
	
	For any $p\in \R^{>0}/\{a\} $ we have 
	$$
	\sigma^{-1}\pi_0^{-1}=\{M\times q^k \ | \ k\in \Z \}.
	$$ 
	and 
	$$
	\sigma_0^{-1} p= \pi \{M\times q^k \ | \ k\in \Z \}\simeq M.
	$$
	Thus, the fibers of $\sigma_0$ are isometric to a statistical manifold $(M,g_M)$ of constant curvature and the manifold $C$ is isomorphic $M\times \R^{>0}/\varphi_q$. 
\end{proof}

	\paragraph*{Acknowledgements.} Many thanks to Misha Verbitsky for fruitful discussions and help with the preparation of the paper.

\end{document}